\documentclass[11pt,reqno]{amsart}
\usepackage{graphicx}
\usepackage{verbatim}
\usepackage{textcomp}
\usepackage{amssymb}
\usepackage{cite}
\usepackage{amsmath}
\usepackage{latexsym}
\usepackage{amscd}
\usepackage{amsthm}
\usepackage{mathrsfs}
\usepackage{bm}
\usepackage{url}
\usepackage{hyperref}
\usepackage{bookmark}
\allowdisplaybreaks[3]
\newtheorem{thm}{Theorem}[section]

\newtheorem{prop}[thm]{Proposition}

\theoremstyle{definition}

\theoremstyle{remark}
\newtheorem{rem}{Remark}[section]
\numberwithin{equation}{section}
\begin{document}
\title[Some integral inequalities]
{Some integral inequalities on weighted Riemannian manifolds with boundary }

\author{Guangyue Huang }

\address{Department of Mathematics,
Henan Normal University, Xinxiang 453007, P.R. China}
\email{hgy@htu.edu.cn(G. Huang) }

\author{Mingfang Zhu }
\email{mfzhu21@126.com(M. Zhu) }

\thanks{The research of authors is supported by NSFC(No. 11971153).}

\begin{abstract}
In this paper, we continue to study some applications with respect to a Reilly type integral formula associated with the $\phi$-Laplacian. Some inequalities of Brascamp-Lieb type and Colesanti type are provided.
\end{abstract}



\maketitle

\section{Introduction}
Let $(M,g)$ be an $n$-dimensional Riemannian manifold with the dimension $n\geq3$, where $g=\langle,\rangle$ is the metric. It is well known that the $\phi$-Laplacian associated with $\phi$ is given by
$$\Delta_\phi v=e^\phi{\rm div}(e^{-\phi}\nabla v)=\Delta v-\langle\nabla\phi,\nabla v\rangle,\ \ \ \ \ \forall \ v\in
C^\infty(M),$$
which is symmetric with respect to the $L^2(M)$ inner product under the weighted measure
$$d\mu=e^{-\phi}dv_g.$$
The $m$-dimensional Bakry-\'{E}mery Ricci curvature (see \cite{Bakry85,li05,Wei09,LW2015,HL2013,HM2016,HZ2016}) associated
with the above $\phi$-Laplacian is given by
\begin{align}\label{Int-1}
{\rm Ric}_{\phi,m}={\rm Ric}+\nabla^2\phi-\frac{1}{m-n}d\phi\otimes d\phi,
\end{align}
where $m$ is a real constant, and $m=n$ if and only if $\phi$ is a constant. Let $m\rightarrow \infty$, then \eqref{Int-1} gives ${\rm Ric}_{\phi,\infty}={\rm Ric}+\nabla^2\phi$. We define
$${\rm Ric}_{\phi}={\rm Ric}_{\phi,\infty}.$$
Thus, ${\rm Ric}_{\phi}$ can be seen as the $\infty$-dimensional Bakry-\'{E}mery Ricci curvature.

For the convenience, we still denote $\Delta,\nabla$ by the Laplacian operator and gradient operator on
$M$, and $\overline{\Delta},\overline{\nabla}$, respectively, by the Laplacian operator and gradient operator on the boundary $\partial M$. The mean curvature $H$ of $\partial M$ is given by
$H={\rm tr}_{g}(II)$, where $II(X,Y)=g(\nabla_X\nu,Y)$ denotes the second fundamental form of $\partial M$
with $\nu$ the outward unit normal on $\partial M$. For any {\em positive} twice differentiable function $V$, we make the following conventions:
\begin{equation}\label{1-Int-1}
{\rm \widehat{Ric}}^{V}_{\phi,m}=\frac{\Delta_\phi V}{V}g-\frac{1}{V}\nabla^2V+{\rm Ric}_{\phi,m},
\end{equation}
$H_\phi=H-\phi_\nu$, $II^V=II-(\ln V)_\nu \overline{g}$ and $d\sigma$ denotes the measure induced on $\partial M$.

In \cite{HMZ2022}, the following Reilly type integral formula associated with $\phi$-Laplacian has been proved:

\vspace*{3mm}
\noindent{\bf Theorem A.} {\em
Let $V$ be a positive twice differentiable function on a given compact Riemannian manifold $M$ with the boundary $\partial M$.
For any smooth function $f$ and $m\in (-\infty,0)\cup [n,+\infty)$, we have
\begin{align}\label{1-th-A-111}
0\leq&\int_{\partial M}\Big[-VII^V\Big(V\overline{\nabla}\frac{f}{V},V\overline{\nabla}\frac{f}{V}\Big)
-V^3H_\phi\Big(\Big(\frac{f}{V}\Big)_\nu\Big)^2\notag\\
&-2V^2\Big(\frac{f}{V}\Big)_\nu\Big(\overline{\Delta}_\phi f-\frac{\overline{\Delta}_\phi V}{V} f\Big)\Big]\,d\sigma+\int_M V\Big[\frac{m-1}{m}\Big(\Delta_\phi f-\frac{\Delta_\phi V}{V}f\Big)^2\notag\\
&-{\rm \widehat{Ric}}^{V}_{\phi,m}\Big(V\nabla\frac{f}{V},V\nabla\frac{f}{V}\Big)\Big]\,d\mu,
\end{align}
where the equality occurs if and only if
\begin{align}\label{Int-4}
\nabla^2 f-\frac{f}{V}\nabla^2 V=\frac{1}{n}\Big(\Delta f-\frac{f}{V}\Delta V\Big)g
\end{align}
and
\begin{align}\label{Int-5}
\Delta f-\frac{f}{V}\Delta V+\frac{n}{\sqrt{(m-n)^2}}V\langle\nabla \phi,\nabla \frac{f}{V}\rangle=0.
\end{align}

}

Using the above theorem, we can obtain the following applications:

\begin{thm}\label{1-thm1}
Let $M$ be a compact Riemannian manifold with a smooth boundary $\partial M$ and $V{\rm \widehat{Ric}}^{V}_{\phi,m}\geq(m-1)Kg$ for some positive constant $K$, where $m\geq n$. If $\lambda$ is the eigenvalue of the problem $\Delta_\phi f-\frac{\Delta_\phi V}{V}f=-\lambda \frac{f}{V}$, then the first nonzero eigenvalue $\lambda_1$ satisfies the following:

(1) In the Dirichlet case, if $H_\phi\geq0$ on $\partial M$, then
$\lambda_1 \geq m K$. Moreover, the equality is attained only when $m=n$, $\phi$ is a constant and $\partial M$ is totally geodesic;

(2) In the Neumann case, if $II^V\geq0$ on $\partial M$, then
$\lambda_1  \geq m K$. Moreover, the equality is attained only when $m=n$ and $\phi$ is a constant.

\end{thm}

\begin{rem}\label{1-rem1}
For a closed Riemannian manifold without boundary, we also obtain the similar results as in Theorem \ref{1-thm1}. On the other hand, we let $L_{\phi}^V(f)=e^{\phi}{\rm div}(e^{-\phi}V^2\nabla\frac{f}{V})$. It is easy to check that the differential operator $L_{\phi}^V$ is linear and the problem $\Delta_\phi f-\frac{\Delta_\phi V}{V}f=-\lambda \frac{f}{V}$ is equivalent to the eigenvalue problem of $L_{\phi}^V(f)=-\lambda f$. In particular, if $u=1$, then Theorem \ref{1-thm1} reduces to the first eigenvalue estimates with respect to $\phi$-Laplacian associated with Bakry-\'{E}mery Ricci curvature.

\end{rem}

Next, we give the following Brascamp-Lieb type inequalities:

\begin{thm}\label{1-thm2}
Let $M$ be a compact Riemannian manifold with a smooth boundary $\partial M$ and ${\rm \widehat{Ric}}^{V}_{\phi,m}>0$, where $m\in (-\infty,0)\cup [n,+\infty)$. Then for any $\varphi\in C^\infty(M)$, we have

(1) In the Dirichlet case, if $H_\phi\geq 0$ and $\varphi=0$ on $\partial M$, then
\begin{align}\label{1-thm2-111}
\frac{m}{m-1}\int_MV\varphi^2\,d\mu\leq\int_MV\langle({\rm \widehat{Ric}}^{V}_{\phi,m})^{-1}\nabla\varphi,\nabla\varphi\rangle\,d\mu.
\end{align}

(2) In the Neumann case, if $II^V\geq0$ on $\partial M$, then
\begin{align}\label{1-thm2-222}
\frac{m}{m-1}\int_MV\Big(\varphi-\frac{\int_MV\varphi\,d\mu}{\int_MV\,d\mu}\Big)^2\,d\mu\leq\int_MV\langle({\rm \widehat{Ric}}^{V}_{\phi,m})^{-1}\nabla\varphi,\nabla\varphi\rangle\,d\mu.
\end{align}

\end{thm}

\begin{rem}
When $V=1$, our Theorem \ref{1-thm2} becomes partial results of Theorem 1.2 of Kolesnikov and Milman in \cite{Kolesnikov2017}.
\end{rem}

\begin{rem}
Taking $\varphi=\frac{f}{V}$ in \eqref{1-thm2-111}, we obtain that if $V{\rm \widehat{Ric}}^{V}_{\phi,m}\geq(m-1)Kg$ for some positive constant $K$, then
\begin{align}\label{2-rem-111}
mK\leq\frac{\int_MV^2|\nabla\frac{f}{V}|^2\,d\mu}{\int_MV(\frac{f}{V})^2\,d\mu}, 
\end{align}
which show that $\lambda_1 \geq m K$ provided $H_\phi\geq 0$ and $f=0$ on $\partial M$. It is exactly the result of (1) of Theorem \ref{1-thm1}.

\end{rem}

As another application, we also give the following Colesanti type inequality:

\begin{thm}\label{1-thm3}
Let $M$ be a compact Riemannian manifold with a smooth boundary $\partial M$ and ${\rm \widehat{Ric}}^{V}_{\phi,m}\geq0$, where $m\in (-\infty,0)\cup [n,+\infty)$. If $II^V>0$ on $\partial M$, then for any $\varphi\in C^\infty(\partial M)$, we have
\begin{align}\label{1-thm3-111}
\int_M V\,d\mu\int_{\partial M}VH_\phi\varphi^2\,d\sigma-&\frac{m-1}{m}\Big(\int_{\partial M}V\varphi\,d\sigma\Big)^2\notag\\
\leq&\int_M V\,d\mu\int_{\partial M}V\langle (II^V)^{-1}\overline{\nabla}\varphi,\overline{\nabla}\varphi\rangle\,d\sigma.
\end{align}

\end{thm}

\begin{rem}
In particular, taking $\varphi=1$ in \eqref{1-thm3-111}, then we obtain
\begin{align}\label{rem-111}
\int_MV\,d\mu\int_{\partial M} V H_\phi \,d\sigma\leq\frac{m-1}{m}\Big(\int_{\partial M}V\,d\sigma\Big)^2,
\end{align}
which is exactly the formula (1.11) of Theorem 1.4 in \cite{HMZ2022}.  On the other hand, using the Cauchy inequality
$$\Big(\int_{\partial M}V\,d\sigma\Big)^2\leq\int_{\partial M}VH_\phi\,d\sigma\int_{\partial M}\frac{V}{H_\phi}\,d\sigma$$
in \eqref{rem-111} gives
\begin{align}\label{rem-222}
\int_MV\,d\mu\leq\frac{m-1}{m}\int_{\partial M}\frac{V}{H_\phi}\,d\sigma,
\end{align}
which is exactly the formula (1.9) of Theorem 1.3 in \cite{HMZ2022}.
\end{rem}

Next, we also can establish the following dual Colesanti type inequality:

\begin{thm}\label{1-thm4}
Let $M$ be a compact Riemannian manifold with a smooth boundary $\partial M$ and $V{\rm \widehat{Ric}}^{V}_{\phi,m}\geq(m-1)Kg$ for some constant $K$, where $m\geq n$. If $H_\phi>0$ on $\partial M$, then for any $\varphi\in C^\infty(\partial M)$, we have
\begin{align}\label{1-thm4-111}
\int_{\partial M}VII^V\Big(V\overline{\nabla}\frac{\varphi}{V},V\overline{\nabla}\frac{\varphi}{V}\Big)\,d\sigma
\leq\int_{\partial M}\frac{1}{H_\phi}V\Big[\Big(\overline{\Delta}_\phi \varphi-\frac{\overline{\Delta}_\phi V}{V} \varphi\Big)+\frac{m-1}{2m}t\frac{\varphi }{V}\Big]^2\,d\sigma,
\end{align}
where $t$ is a parameter satisfying $t\leq mK$.
\end{thm}

\begin{rem}
Obviously, our Theorem \ref{1-thm4} generalizes Theorem 1.2 of Kolesnikov and Milman in \cite{Kolesnikov2018}.
\end{rem}

\section{Proof of results}
\subsection{Proof of Theorem \ref{1-thm1}}
Using the divergence theorem, we have
\begin{align}\label{1-Sec-Prof-Th-1}
&\int_M V\Big(\Delta_\phi f-\frac{\Delta_\phi V}{V}f\Big)^2\,d\mu\notag\\
=&\int_M\Big(\Delta_\phi f-\frac{\Delta_\phi V}{V}f\Big)(V\Delta_\phi f-f\Delta_\phi V)\,d\mu\notag\\
=&\int_M\Big(\Delta_\phi f-\frac{\Delta_\phi V}{V}f\Big)\Big[e^{-\phi}V^2\Big(\frac{f}{V}\Big)_i\Big]_i\,dv_g\notag\\
=&\int_{\partial M}V^2\Big(\Delta_\phi f-\frac{\Delta_\phi V}{V}f\Big)\Big(\frac{f}{V}\Big)_\nu\,d\sigma-\int_MV^2\Big\langle \nabla\frac{f}{V},\nabla\Big(\Delta_\phi f-\frac{\Delta_\phi V}{V}f\Big)\Big\rangle\,d\mu.
\end{align}
On the other hand, by virtue of
$$\Delta_\phi f=\overline{\Delta}_\phi f+H_\phi f_{\nu}+f_{\nu\nu}$$
and
$$\Delta_\phi V=\overline{\Delta}_\phi V+H_\phi V_{\nu}+V_{\nu\nu},$$
we obtain that on $\partial M$, it holds that
\begin{align}\label{1-Sec-Prof-Th-2}
\Delta_\phi f-\frac{\Delta_\phi V}{V}f=&\overline{\Delta}_\phi f-\frac{\overline{\Delta}_\phi V}{V} f+V H_\phi \Big(\frac{f}{V}\Big)_\nu+f_{\nu\nu}-\frac{f}{V}V_{\nu\nu}.
\end{align}
This shows that \eqref{1-Sec-Prof-Th-1} can be written as
\begin{align}\label{1-Sec-Prof-Th-3}
&\int_M V\Big(\Delta_\phi f-\frac{\Delta_\phi V}{V}f\Big)^2\,d\mu\notag\\
=&\int_{\partial M}\Big[V^2\Big(\frac{f}{V}\Big)_\nu\Big(\overline{\Delta}_\phi f-\frac{\overline{\Delta}_\phi V}{V} f\Big)+V^3 H_\phi\Big(\Big(\frac{f}{V}\Big)_\nu\Big)^2\notag\\
&+V^2\Big(f_{\nu\nu}-\frac{f}{V}V_{\nu\nu}\Big)\Big(\frac{f}{V}\Big)_\nu\Big]\,d\sigma-\int_MV^2\Big\langle\nabla\frac{f}{V},\nabla\Big(\Delta_\phi f-\frac{\Delta_\phi V}{V}f\Big)\Big\rangle\,d\mu.
\end{align}
Hence, putting \eqref{1-Sec-Prof-Th-3} into \eqref{1-th-A-111} gives an alternative form of Reilly type integral formula:

\begin{prop}\label{1-prop1}
Let $V$ be a positive twice differentiable function on a given compact Riemannian manifold $M$ with the boundary $\partial M$.
For any smooth function $f$ and $m\in (-\infty,0)\cup [n,+\infty)$, we have
\begin{align}\label{1-prop-111}
0\leq&\int_{\partial M}\Big[-VII^V\Big(V\overline{\nabla}\frac{f}{V},V\overline{\nabla}\frac{f}{V}\Big)
-\frac{1}{m}V^3H_\phi\Big(\Big(\frac{f}{V}\Big)_\nu\Big)^2\notag\\
&-\frac{m+1}{m}V^2\Big(\frac{f}{V}\Big)_\nu\Big(\overline{\Delta}_\phi f-\frac{\overline{\Delta}_\phi V}{V} f\Big)+\frac{m-1}{m}V^2\Big(f_{\nu\nu}-\frac{f}{V}V_{\nu\nu}\Big)\Big(\frac{f}{V}\Big)_\nu\Big]\,d\sigma\notag\\
&-\int_M \Big[\frac{m-1}{m}V^2\Big\langle\nabla\frac{f}{V},\nabla\Big(\Delta_\phi f-\frac{\Delta_\phi V}{V}f\Big)\Big\rangle+V{\rm \widehat{Ric}}^{V}_{\phi,m}\Big(V\nabla\frac{f}{V},V\nabla\frac{f}{V}\Big)\Big]\,d\mu,
\end{align}
where the equality occurs if and only if \eqref{Int-4} and \eqref{Int-5} occur.

\end{prop}

Now, we are in the position to complete the proof of Theorem \ref{1-thm1}.
In the Dirichlet case, we let $f$ is a solution to
$$\Delta_\phi f-\frac{\Delta_\phi V}{V}f=-\lambda  \frac{f}{V}\  {\rm on }\ M,\ \ \ f=0\ {\rm on }\ \partial M.$$
Since the above problem is equivalent to $e^{\phi}{\rm div}(e^{-\phi}V^2\nabla\frac{f}{V})=-\lambda f$, the existence and uniqueness is due to the Fredholm alternative (for detail, see Remark \ref{1-rem1} or page 512 in \cite{lixia2019}). Hence, from \eqref{1-prop-111}, we obtain
\begin{align}\label{1-Sec-Prof-Th-4}
0\leq&\int_{\partial M}\Big[-\frac{1}{m}V^3H_\phi\Big(\Big(\frac{f}{V}\Big)_\nu\Big)^2+\frac{m-1}{m}V^2\Big(f_{\nu\nu}
-\frac{f}{V}V_{\nu\nu}\Big)\Big(\frac{f}{V}\Big)_\nu\Big]\,d\sigma\notag\\
&-\int_M \Big[\frac{m-1}{m}V^2\Big\langle\nabla\frac{f}{V},\nabla\Big(\Delta_\phi f-\frac{\Delta_\phi V}{V}f\Big)\Big\rangle+V{\rm \widehat{Ric}}^{V}_{\phi,m}(V\nabla\frac{f}{V},V\nabla\frac{f}{V})\Big]\,d\mu\notag\\
\leq&\int_{\partial M}\Big[-\frac{1}{m}V^3H_\phi\Big(\Big(\frac{f}{V}\Big)_\nu\Big)^2+\frac{m-1}{m}V^2\Big(f_{\nu\nu}
-\frac{f}{V}V_{\nu\nu}\Big)\Big(\frac{f}{V}\Big)_\nu\Big]\,d\sigma\notag\\
&+\frac{m-1}{m}(\lambda-mK)\int_MV^2\Big|\nabla\frac{f}{V}\Big|^2\,d\mu.
\end{align}
Since $f|_{\partial M}=0$, it follows from \eqref{1-Sec-Prof-Th-2} that
$f_{\nu\nu}-\frac{f}{V}V_{\nu\nu}=-V H_\phi \Big(\frac{f}{V}\Big)_\nu$. Then \eqref{1-Sec-Prof-Th-4} becomes
\begin{align}\label{1-Sec-Prof-Th-5}
0\leq&-\int_{\partial M}V^3H_\phi\Big(\Big(\frac{f}{V}\Big)_\nu\Big)^2\,d\sigma+\frac{m-1}{m}(\lambda-mK)\int_MV^2\Big|\nabla\frac{f}{V}\Big|^2\,d\mu\notag\\
\leq&\frac{m-1}{m}(\lambda-mK)\int_MV^2\Big|\nabla\frac{f}{V}\Big|^2\,d\mu
\end{align}
from $H_\phi\geq0$, which shows that $\lambda\geq m K g$. When $\lambda= m K$, we have $V{\rm \widehat{Ric}}^{V}_{\phi,m}=(m-1)K$ and all inequalities must be equal. If $m>n$, then \eqref{Int-5} shows that
$$V\Delta f-f\Delta V+\frac{n}{m-n}V^2\langle\nabla \phi,\nabla \frac{f}{V}\rangle=0,$$
which is equivalent to
\begin{align}\label{1-Sec-Prof-Th-6}
e^{-\frac{n}{m-n}\phi}\Big[e^{\frac{n}{m-n}\phi}V^2\Big(\frac{f}{V}\Big)_i\Big]_i=&\Big[V^2\Big(\frac{f}{V}\Big)_i\Big]_i
+\frac{n}{m-n}V^2\langle\nabla \phi,\nabla \frac{f}{V}\rangle\notag\\
=&0,
\end{align}
where we notice that $V\Delta \frac{f}{V}=\Delta f-\frac{f}{V}\Delta V-2\langle\nabla V,\nabla \frac{f}{V}\rangle$.
Multiplying both sides of \eqref{1-Sec-Prof-Th-6} with $\frac{f}{V}e^{\frac{n}{m-n}\phi}$ gives
\begin{align}\label{1-Sec-Prof-Th-7}
0=&\int_M\frac{f}{V}\Big[e^{\frac{n}{m-n}\phi}V^2\Big(\frac{f}{V}\Big)_i\Big]_i\,dv_g\notag\\
=&-\int_M V^2\Big|\nabla\frac{f}{V}\Big|^2 e^{\frac{n}{m-n}\phi}\,dv_g,
\end{align}
and then $f=\theta V$, where $\theta$ is a constant. This contradicts with that $f$ is nontrivial, which shows that $m=n$ and $\phi$ must be constant. In this case, \eqref{Int-4} implies
\begin{align}\label{1-Sec-Prof-Th-8}
\nabla^2 f-\frac{f}{V}\nabla^2 V=-K\frac{f}{V}g.
\end{align}
Restrict \eqref{1-Sec-Prof-Th-8} on $\partial M$ and using $f|_{\partial M}=0$, we have $II_{\alpha\beta}f_{\nu}=0$ and then $\partial M$ is totally geodesic.

In the Neumann case, we let $f$ is a solution to
$$\Delta_\phi f-\frac{\Delta_\phi V}{V}f=-\lambda  \frac{f}{V}\  {\rm on }\ M,\ \ \ V\Big(\frac{f}{V}\Big)_\nu=0\ {\rm on }\ \partial M.$$
Then, from \eqref{1-prop-111}, we obtain
\begin{align}\label{1-Sec-Prof-Th-10}
0\leq&-\int_{\partial M}VII^V\Big(V\overline{\nabla}\frac{f}{V},V\overline{\nabla}\frac{f}{V}\Big)
\,d\sigma\notag\\
&-\int_M \Big[\frac{m-1}{m}V^2\Big\langle\nabla\frac{f}{V},\nabla\Big(\Delta_\phi f-\frac{\Delta_\phi V}{V}f\Big)\Big\rangle+V{\rm \widehat{Ric}}^{V}_{\phi,m}\Big(V\nabla\frac{f}{V},V\nabla\frac{f}{V}\Big)\Big]\,d\mu\notag\\
\leq&\frac{m-1}{m}(\lambda-mK)\int_MV^2\Big|\nabla\frac{f}{V}\Big|^2\,d\mu
\end{align}
from $II^V\geq0$, which shows that $\lambda\geq m K$. Similarly, we can also prove that $m=n$ and $\phi$ must be constant.

\subsection{Proof of Theorem \ref{1-thm2}}
(1) In the Dirichlet case. We assume that $\varphi=0$ on $\partial M$, and solve the following Dirichlet poisson equation:
$$\Delta_\phi f-\frac{\Delta_\phi V}{V}f=\varphi \  {\rm on }\ M,\ \ \ f=0\ {\rm on }\ \partial M.$$
Thus, it follows from \eqref{1-th-A-111} that
\begin{align}\label{2-thm-Proof-1}
\frac{m-1}{m}\int_MV\Big(\Delta_\phi f-\frac{\Delta_\phi V}{V}f\Big)^2\,d\mu\geq\int_MV^3\langle{\rm \widehat{Ric}}^{V}_{\phi,m}(\nabla\frac{f}{V}),\nabla\frac{f}{V}\rangle\,d\mu,
\end{align}
which is equivalent to
\begin{align}\label{2-thm-Proof-2}
\frac{m-1}{m}\int_MV\varphi^2\,d\mu\geq\int_MV^3\langle{\rm \widehat{Ric}}^{V}_{\phi,m}(\nabla\frac{f}{V}),\nabla\frac{f}{V}\rangle\,d\mu.
\end{align}
Using the divergence theorem, we have
\begin{align}\label{2-thm-Proof-3}
\int_MV\varphi^2\,d\mu=&\int_M\varphi(V\Delta_\phi f-f\Delta_\phi V)\,d\mu\notag\\
=&-\int_MV^2\langle\nabla\varphi,\nabla\frac{f}{V}\rangle\,d\mu\notag\\
\leq&\Big[\int_MV^3\langle{\rm \widehat{Ric}}^{V}_{\phi,m}(\nabla\frac{f}{V}),\nabla\frac{f}{V}\rangle\,d\mu\Big]^{\frac{1}{2}}\Big[\int_MV\langle({\rm \widehat{Ric}}^{V}_{\phi,m})^{-1}\nabla\varphi,\nabla\varphi\rangle\,d\mu\Big]^{\frac{1}{2}}.
\end{align}
Thus, putting \eqref{2-thm-Proof-2} into \eqref{2-thm-Proof-3} gives
the desired estimate \eqref{1-thm2-111}.

(2) In the Neumann case, for any $\varphi\in C^\infty(M)$ with $\int_MV\varphi\,d\mu=0$, there exists a function $f$ such that
$$\Delta_\phi f-\frac{\Delta_\phi V}{V}f=\varphi \  {\rm on }\ M,\ \ \ V\Big(\frac{f}{V}\Big)_\nu=0\ {\rm on }\ \partial M.$$
Thus, by using \eqref{1-th-A-111}, we also have that \eqref{2-thm-Proof-2} holds. Using the Neumann condition, we also have that \eqref{2-thm-Proof-3} is true and the estimate
\begin{align}\label{2-thm-Proof-4}
\frac{m}{m-1}\int_MV\varphi^2\,d\mu\leq\int_MV\langle({\rm \widehat{Ric}}^{V}_{\phi,m})^{-1}\nabla\varphi,\nabla\varphi\rangle\,d\mu
\end{align}
follows by insetting \eqref{2-thm-Proof-2} into \eqref{2-thm-Proof-3}.
In particular, we let
$$\widetilde{\varphi}=\varphi-\frac{\int_MV\varphi\,d\mu}{\int_MV\,d\mu}.$$
Then we have $\int_MV\widetilde{\varphi}\,d\mu=0$ and the estimate \eqref{1-thm2-222} follows.

\subsection{Proof of Theorem \ref{1-thm3}}
We consider the following Neumann problem:
$$\Delta_\phi f-\frac{\Delta_\phi V}{V}f=C \  {\rm on }\ M,\ \ \ V\Big(\frac{f}{V}\Big)_\nu=\varphi\ {\rm on }\ \partial M,$$
which implies that
$$C=\frac{\int_{\partial M}V\varphi\,d\sigma}{\int_MV\,d\mu}.$$
Then, from \eqref{1-th-A-111} we obtain
\begin{align}\label{3-thm-Proof-1}
\frac{m-1}{m}C^2\int_M V\,d\mu\geq&\int_{\partial M}\Big[VII^V\Big(V\overline{\nabla}\frac{f}{V},V\overline{\nabla}\frac{f}{V}\Big)
+VH_\phi\varphi^2\notag\\
&+2V\varphi\Big(\overline{\Delta}_\phi f-\frac{\overline{\Delta}_\phi V}{V} f\Big)\Big]\,d\sigma\notag\\
=&\int_{\partial M}\Big[V^3II^V\Big(\overline{\nabla}\frac{f}{V},\overline{\nabla}\frac{f}{V}\Big)
+VH_\phi\varphi^2\notag\\
&-2V^2\langle\overline{\nabla}\varphi,\overline{\nabla}\frac{f}{V}\rangle\Big]\,d\sigma.
\end{align}
Applying the inequality
$$-2V^2\langle\overline{\nabla}\varphi,\overline{\nabla}\frac{f}{V}\rangle\geq-V^3\langle II^V(\overline{\nabla}\frac{f}{V}),\overline{\nabla}\frac{f}{V}\rangle-V\langle (II^V)^{-1}\overline{\nabla}\varphi,\overline{\nabla}\varphi\rangle$$
in \eqref{3-thm-Proof-1} yields
\begin{align}\label{3-thm-Proof-2}
\frac{m-1}{m}C^2\int_M V\,d\mu\geq&\int_{\partial M}\Big[VH_\phi\varphi^2-V\langle (II^V)^{-1}\overline{\nabla}\varphi,\overline{\nabla}\varphi\rangle\Big]\,d\sigma,
\end{align}
which is equivalent to \eqref{1-thm3-111}.

\subsection{Proof of Theorem \ref{1-thm4}}
We consider the following Dirichlet problem:
$$\Delta_\phi f-\frac{\Delta_\phi V}{V}f+t \frac{f}{V}=0 \  {\rm on }\ M,\ \ \ f=\varphi\ {\rm on }\ \partial M,$$
where $t\leq mK$.
Then \eqref{1-prop-111} shows that
\begin{align}\label{4-thm-Proof-1}
0\geq&\int_{\partial M}\Big[VII^V\Big(V\overline{\nabla}\frac{f}{V},V\overline{\nabla}\frac{f}{V}\Big)
+\frac{1}{m}V^3H_\phi\Big(\Big(\frac{f}{V}\Big)_\nu\Big)^2\notag\\
&+\frac{m+1}{m}V^2\Big(\frac{f}{V}\Big)_\nu\Big(\overline{\Delta}_\phi f-\frac{\overline{\Delta}_\phi V}{V} f\Big)-\frac{m-1}{m}V^2\Big(f_{\nu\nu}-\frac{f}{V}V_{\nu\nu}\Big)\Big(\frac{f}{V}\Big)_\nu\Big]\,d\sigma\notag\\
&+\frac{m-1}{m}(mK-t)\int_M V^2\Big|\nabla\frac{f}{V}\Big|^2\,d\mu\notag\\
\geq&\int_{\partial M}\Big[VII^V\Big(V\overline{\nabla}\frac{f}{V},V\overline{\nabla}\frac{f}{V}\Big)
+\frac{1}{m}V^3H_\phi\Big(\Big(\frac{f}{V}\Big)_\nu\Big)^2\notag\\
&+\frac{m+1}{m}V^2\Big(\frac{f}{V}\Big)_\nu\Big(\overline{\Delta}_\phi f-\frac{\overline{\Delta}_\phi V}{V} f\Big)-\frac{m-1}{m}V^2\Big(f_{\nu\nu}-\frac{f}{V}V_{\nu\nu}\Big)\Big(\frac{f}{V}\Big)_\nu\Big]\,d\sigma.
\end{align}
By virtue of \eqref{1-Sec-Prof-Th-2}, we have
\begin{align}\label{4-thm-Proof-2}
-t\frac{f}{V}=&\overline{\Delta}_\phi f-\frac{\overline{\Delta}_\phi V}{V} f+V H_\phi \Big(\frac{f}{V}\Big)_\nu+f_{\nu\nu}-\frac{f}{V}V_{\nu\nu},
\end{align}
and then \eqref{4-thm-Proof-1} becomes
\begin{align}\label{4-thm-Proof-3}
0\geq&\int_{\partial M}\Big[VII^V\Big(V\overline{\nabla}\frac{f}{V},V\overline{\nabla}\frac{f}{V}\Big)
+V^3H_\phi\Big(\Big(\frac{f}{V}\Big)_\nu\Big)^2\notag\\
&+2V^2\Big(\frac{f}{V}\Big)_\nu\Big(\overline{\Delta}_\phi f-\frac{\overline{\Delta}_\phi V}{V} f\Big)+\frac{m-1}{m}t\varphi V\Big(\frac{f}{V}\Big)_\nu\Big]\,d\sigma\notag\\
=&\int_{\partial M}\Big[VII^V\Big(V\overline{\nabla}\frac{\varphi}{V},V\overline{\nabla}\frac{\varphi}{V}\Big)
+V^3H_\phi\Big(\Big(\frac{f}{V}\Big)_\nu\Big)^2\notag\\
&+V^2\Big(2\Big(\overline{\Delta}_\phi \varphi-\frac{\overline{\Delta}_\phi V}{V} \varphi\Big)+\frac{m-1}{m}t\frac{\varphi }{V}\Big)\Big(\frac{f}{V}\Big)_\nu\Big]\,d\sigma.
\end{align}
Using the inequality $ax^2+bx\geq-\frac{b^2}{4a}$ with $a>0$, we obtain
\begin{align}\label{4-thm-Proof-4}
V^3H_\phi\Big(\Big(\frac{f}{V}\Big)_\nu\Big)^2+&V^2\Big(2\Big(\overline{\Delta}_\phi \varphi-\frac{\overline{\Delta}_\phi V}{V} \varphi\Big)+\frac{m-1}{m}t\frac{\varphi }{V}\Big)\Big(\frac{f}{V}\Big)_\nu\notag\\
\geq&-\frac{1}{H_\phi}V\Big[\Big(\overline{\Delta}_\phi \varphi-\frac{\overline{\Delta}_\phi V}{V} \varphi\Big)+\frac{m-1}{2m}t\frac{\varphi }{V}\Big]^2.
\end{align}
Thus, putting \eqref{4-thm-Proof-4} into \eqref{4-thm-Proof-3} completes the proof of Theorem \ref{1-thm4}.

\section{Further remarks}
Let $\Omega^n\subset \mathbb{H}^n$ ($\mathbb{S}_+^n$, respectively) be a compact domain with smooth boundary $\partial \Omega$. Let $V=\cosh r$, $K=-1$ or $V=\cos r$, $K=1$ for the case $H^n$ or $S_+^n$, respectively, where $r(x)={\rm dist} (x, p)$ is the distance function for the fixed point $p\in\mathbb{H}^n$ ($p\in\mathbb{S}_+^n$, respectively), see \cite{Xia2016}. Then we have
\begin{align}\label{Sec-3-Proof-1}
\nabla^2V=-KVg.
\end{align}
In this case, we obtain $V{\rm \widehat{Ric}}^{V}:=V{\rm \widehat{Ric}}^{V}_{0,n}=(\Delta V) g-\nabla^2V+{\rm Ric}=0$. Hence, Theorems \ref{1-thm3} and \ref{1-thm4} can be stated as follows:

\begin{thm}\label{Rem-thm1}
Let $\Omega\subset \mathbb{H}^n$ ($\mathbb{S}_+^n$, respectively) be a compact domain with smooth boundary $\partial \Omega$ and $V$ defined as above. If $II^V>0$ on $\partial \Omega$, then for any $\varphi\in C^\infty(\partial \Omega)$, we have
\begin{align}\label{Rem-thm1-thm3-111}
\int_\Omega V\,dv_g\int_{\partial \Omega}VH \varphi^2\,dv_{\overline{g}}-&\frac{n-1}{n}\Big(\int_{\partial \Omega}V\varphi\,dv_{\overline{g}}\Big)^2\notag\\
\leq&\int_\Omega V\,dv_g\int_{\partial \Omega}V\langle (II^V)^{-1}\overline{\nabla}\varphi,\overline{\nabla}\varphi\rangle\,dv_{\overline{g}}.
\end{align}

\end{thm}

\begin{thm}\label{Rem-thm4}
Let $\Omega\subset \mathbb{H}^n$ ($\mathbb{S}_+^n$, respectively) be a compact domain with smooth boundary $\partial \Omega$ and $V$ defined as above. If $H>0$ on $\partial \Omega$, then for any $\varphi\in C^\infty(\partial \Omega)$, we have
\begin{align}\label{Rem-thm1-thm4-111}
\int_{\partial \Omega}VII^V\Big(V\overline{\nabla}\frac{\varphi}{V},V\overline{\nabla}\frac{\varphi}{V}\Big)\,dv_{\overline{g}}
\leq\int_{\partial \Omega}\frac{1}{H}V\Big(\overline{\Delta} \varphi-\frac{\overline{\Delta} V}{V} \varphi\Big)^2\,dv_{\overline{g}}.
\end{align}

\end{thm}

In particular, taking $\varphi=1$ in \eqref{Rem-thm1-thm3-111}, we have
\begin{align}\label{Sec-3-Proof-2}
\int_\Omega V\,dv_g\int_{\partial \Omega} V H \,dv_{\overline{g}}\leq\frac{n-1}{n}\Big(\int_{\partial \Omega}V\,dv_{\overline{g}}\Big)^2,
\end{align}
which is exactly the formula (5) of Theorem 1.1 in \cite{Xia2016}. Therefore, our
Theorem \ref{Rem-thm1} generalizes Theorem 1.1 of Xia in \cite{Xia2016}.

On the other hand, for a closed hypersurface $x:(\Sigma, g)\rightarrow \mathbb{H}^{n+1}$ ($\mathbb{S}_+^{n+1}$, respectively) and $V$ defined as above, the second fundamental form of $x$ defined by
\begin{align}\label{Sec-3-Proof-3}
II_x^V=II_x-(\ln V)_\nu g
\end{align}
satisfies $II_x^V\geq0$. Moreover, in the case of $\mathbb{H}^{n+1}$ we assume that $II_x^V>0$ and $\Sigma$ is horo-convex and in the case of $\mathbb{S}_+^{n+1}$ we assume that $\Sigma$ is convex, then we have
\begin{align}\label{Sec-3-Proof-4}
V{\rm \widehat{Ric}}^{V}_{x}:=(\overline{\Delta} V)g -\overline{\nabla}V+V{\rm Ric}_{\Sigma}>0,
\end{align}
where ${\rm Ric}_{\Sigma}$ is the Ricci curvature of $\Sigma$ (see page 515 in \cite{lixia2019}).

Hence, applying Theorem \ref{1-thm2} on $\Sigma$ yields the following:

\begin{thm}\label{Rem-thm2}
Let $x:(\Sigma, g)\rightarrow \mathbb{H}^{n+1}$ ($\mathbb{S}_+^{n+1}$, respectively) be a closed hypersurface.
In the case of $\mathbb{H}^{n+1}$ we assume that $II_x^V>0$ and $\Sigma$ is horo-convex. In the case of $\mathbb{S}_+^{n+1}$ we assume that $\Sigma$ is convex. Then for any $\varphi\in C^\infty(\Sigma)$, the following inequality holds:
\begin{align}\label{Rem-thm1-thm4-111}
\frac{n}{n-1}\int_\Sigma V\Big(\varphi-\frac{\int_MV\varphi\,d\mu}{\int_MV\,d\mu}\Big)^2\,dv_g\leq\int_\Sigma V\langle({\rm \widehat{Ric}}^{V}_{x})^{-1}\nabla\varphi,\nabla\varphi\rangle\,dv_g.
\end{align}

\end{thm}

\bibliographystyle{Plain}

\begin{thebibliography}{10}

\bibitem{Bakry85}
D. Bakry, M. Emery,
Diffusion hypercontractives, S\'{e}m. Prob. XIX. Lect. Notes in Math.,
vol. 1123, pp. 177-206 (1985).


\bibitem{HMZ2022}
G.Y. Huang, B.Q. Ma, M.F. Zhu,
A Reilly type integral formula and its applications, arXiv:2201.09439

\bibitem{HL2013}
G.Y. Huang, H.Z, Li,
Gradient estimates and entropy formulae of porous medium and fast diffusion equations for the Witten Laplacian,
Pacific J. Math., 268(2014), 47-78.

\bibitem{HZ2016}
G.Y. Huang, F.Q, Zeng,
De Lellis-Topping type inequalities for $f$-Laplacians,
Studia Math., 232(2016), 189-199.

\bibitem{HM2016}
G.Y. Huang, B.Q, Ma,
Sharp bounds for the first nonzero Steklov eigenvalues for $f$-Laplacians,
Turkish J. Math., 40(2016), 770-783.


\bibitem{Kolesnikov2017}
A.V. Kolesnikov, E. Milman,
Brascamp-Lieb-type inequalities on weighted Riemannian manifolds with boundary,
J. Geom. Anal., 27(2017), 1680-1702.

\bibitem{Kolesnikov2018}
A.V. Kolesnikov, E. Milman,
Poincar\'{e} and Brunn-Minkowski inequalities on the boundary of weighted Riemannian manifolds,
Amer. J. Math., 140(2018), 1147-1185.

\bibitem{li05}
X.-D. Li,
Liouville theorems for symmetric diffusion operators on complete Riemannian manifolds,
J. Math. Pures Appl., 84(2005), 1361-1995.

\bibitem{LW2015}
H.Z, Li, Y. Wei,
$f$-minimal surface and manifold with positive $m$-Bakry-\'{E}mery Ricci curvature,
J. Geom. Anal., 25(2015), 421-435.

\bibitem{lixia2019}
J.F. Li, C. Xia,
An integral formula and its applications on sub-static manifolds,
J. Differential Geom., 113(2019), 493-518.

\bibitem{MW2016}
P.Z. Miao, X.D. Wang,
Boundary effect of Ricci curvature,
J. Differential Geom., 103(2016), 59-82.

\bibitem{Xia2016}
C. Xia,
A Minkowski type inequality in space forms,
Calc. Var. Partial Differential Equations, 55(2016),  Art. 96, 8 pp.


\bibitem{Wei09}
G.-F. Wei, W. Wylie,
Comparison geometry for the Bakry-Emery Ricci tensor,
J. Differential Geom., 83(2009), 377-405.



\end{thebibliography}

\end{document}